\documentclass{amsart}
\usepackage{amssymb}
\usepackage{tikz}
\usepackage{color}
\usepackage[mathscr]{eucal}
\usepackage{hyperref}
\usepackage{wasysym}
\usepackage{tikz}

\newtheorem{theorem}{Theorem}[section]

\theoremstyle{definition}

\numberwithin{equation}{section}

\def\ggg{\mathfrak{g}}

\def\ggg{\mathfrak{g}}

\def\hhh{\mathfrak{h}}

\def\bbc{\mathbb{C}}

\def\bbz{\mathbb{Z}}

\def\bo{{\bar 1}}
\def\bz{{\bar 0}}

\def\wsc{\widehat{\mathscr{C}}}

\def\spo{\textsf{S}_{\Pi_0}}
\def\scrw{\mathscr{W}}

\begin{document}

\title{Coxeter graphs for super weyl groups of exceptional classical Lie superalgebras}\label{se:6}

\author{Yuhui shen  and Zhiyang Tan}

\address{School of Mathematical Sciences, Ministry of Education Key Laboratory of Mathematics and Engineering Applications \& Shanghai Key Laboratory of PMMP,  East China Normal University, No. 500 Dongchuan Rd., Shanghai 200241, China}\email{52275500050@stu.ecnu.edu.cn}
\address{School of Mathematical Sciences, Ministry of Education Key Laboratory of Mathematics and Engineering Applications \& Shanghai Key Laboratory of PMMP,  East China Normal University, No. 500 Dongchuan Rd., Shanghai 200241, China}\email{51275500002@stu.ecnu.edu.cn}

\begin{abstract}
	Super Weyl group plays an important role in the study of representations of basic classical Lie superalgebras. The Coxeter graphs for super Weyl groups of basis classical Lie superalgebras have been given in \cite{CLS}, where the authors also made a proposal on the Coxeter graphs for the super Weyl groups of exceptional classical Lie superalgebras $D(2,1,\alpha)$, $F(4)$ and $G(3)$. In this paper, we present all fundamental systems of the exceptional Lie superalgebras and get the Coxeter graphs for the corresponding super Weyl groups, verifying a proposal in \cite{CLS}.
	\end{abstract}
\maketitle

\section{Introduction}
It is well-known that all Cartan subalgebras of a semisimple Lie algebra over the complex number field are conjugate and all Borel subalgebras are also conjugate. The Weyl group plays an important role in the representation theory of complex semisimple Lie algebras. However, it is not the case for Lie superalgebras. The Weyl groups for basical classical Lie superalgebras only reflect partial information in representation theory of those Lie superalgebras. In \cite{PS}, the authors introduced the super Weyl group $\wsc$ for a basical Lie superalgebra, 
which is the transformation group of the set of all Borel subalgebras containing a given Cartan subalgebra,  generated by ordinary reflections and odd reflections. The introduction of super Weyl groups integrates ordinary reflections and odd reflections in a body. With aid of super Weyl groups,  the authors in \cite{PS} and \cite{LS} successfully  obtained the results on Jantzen filtrations and Jantzen sum formulas in  modular representations of basic classical Lie algebras and general linear supergroups, respectively.   This demonstrates that super Weyl groups deserve further study.  Generally,  a super Weyl group must be a quotient of the corresponding Coxeter group.  In \cite{CLS}, the authors initiate to study the Coxeter systems for super Weyl groups, by drawing the Coxeter graphs for basic classical Lie superalgebras of type $A$, $B$, $C$ and $D$. However, they didn't finish the arguments for the exceptional cases.

In this paper, we continue to study the Coxeter system  for super Weyl groups of exceptional cases. We also follow the notaions in \cite{CLS} and denote super Weyl group by $\wsc$. For super case, the super Weyl groups are independent of the choices of Borel subalgebras and Weyl groups can be realized as the subgroups of super Weyl groups. The purpose of this paper is to calculate the $m_{xy}$ defined in \cite{CLS}, where $x$ is an odd reflection and get the Coxeter graphs for exceptional classical Lie superalgebras $D(2,1,\alpha)$, $F(4)$ and $G(3)$. In order to calculate $m_{xy}$, we need to first describe all fundamental systems and consider the action of the multiplication of odd reflection and ordinary reflection on the fundamental systems. There are too many fundamental systems $\widetilde\Pi^{'}$ for the exceptional classical Lie superalgebras and the computation will be very tedious. However, one can see that if a fundamental system $\widetilde\Pi^{'}$ contains some simple root $\alpha$ transformed into the corresponding isotropic odd root up to a sign under the reflections, then the order on this fundamental system is bigger than $2$. Otherwise, the order will be $2$. Fortunately, we can classify the fundamental systems into some cases to simplify the calculation of $m_{xy}$. Finally, the Coxeter groups for $D(2,1,\alpha)$, $F(4)$ and $G(3)$ can be obtained.

In this paper, we directly begin with our arguments on Coxeter graphs for super Weyl groups. The reader can refer to \cite{CW}, \cite{K} for exceptional super Lie superalgebras,  to \cite{CSW}, \cite{CW}, or \cite{K} for their root systems, and to  \cite{GP} or \cite{Hum1, Hum2} for Coxeter groups.

\section{The root system and Coxeter graph for $D(2,1,\alpha)$} \subsection{The root system}


Recall that the Lie superalgebra $\ggg = D(2|1; \alpha)$ is a family of simple contragredient Lie superalgebras, which depends on a parameter $\alpha$.
There are isomorphisms of Lie superalgebras
with different parameters
$D(2|1; \alpha)\cong D(2,1;-1 -\alpha^{-1})\cong D(2|1; \alpha^{-1})$ when $\alpha$ is nonzero.

Let $\ggg=D(2,1;\alpha)$. Recall that $\ggg= \ggg_\bz\oplus \ggg_\bo$ with $\ggg_\bz\cong \mathfrak{sl}(2)^{\oplus 3}$ and as a $\ggg_\bz$-module, $\ggg_\bo\cong (\bbc^2)^{\oplus 2}$.
For a Cartan subalgebra $\hhh$ of $\ggg_\bz$, its linear dual $\hhh^*$ has a basis $\{\delta, \varepsilon_{1}, \varepsilon_{2} \}$ with $(\delta,\delta)=-(1+\alpha), (\varepsilon_{1},\varepsilon_{1})=1, (\varepsilon_{2},\varepsilon_{2})=\alpha, (\varepsilon_{1},\varepsilon_{2})=0$ and $(\delta,\varepsilon_{i})=0, i=1,2$ (see for example, \cite{K} or \cite{CSW} for the details).
Its root system
$\Phi=\Phi_{\bar{0}}\bigcup \Phi_{\bar{1}}$ is presented as
$$\{\pm2\delta,\pm2\varepsilon_{i}\}\cup
\{\pm\delta\pm\varepsilon_{1}\pm\varepsilon_{2}\}.$$
    We can take
     $$\Pi=\{\delta+\varepsilon_1+ \varepsilon_2, -2\varepsilon_1, -2\varepsilon_2\}$$
as the standard fundamental system of $D(2,1,\alpha)$.  Then the corresponding extended standard fundamental system is
   $$\widetilde \Pi=\{-2\delta\}\cup \{ \delta+\varepsilon_1+ \varepsilon_2, -2\varepsilon_1, -2\varepsilon_2\}.$$
What's more, we can classify all fundamental systems $\widetilde\Pi^{'}$ of $D(2,1,\alpha)$ into following $4$ cases:

\begin{align*}
& \widetilde\Pi^{'}=\{s\delta+t_{1}\varepsilon_{1}+t_{2}\varepsilon_{2},-2t_{1}\varepsilon_{1},-2t_{2}\varepsilon_{2}\} \cr
& \widetilde\Pi^{'}=\{s\delta+t_{1}\varepsilon_{1}+t_{2}\varepsilon_{2},-s\delta+t_{1}\varepsilon_{1}-t_{2}\varepsilon_{2},-s\delta-t_{1}\varepsilon_{1}+t_{2}\varepsilon_{2}\}\cr
&
\widetilde\Pi^{'}=\{2t_{1}\varepsilon_{1},s\delta-t_{1}\varepsilon_{1}+t_{2}\varepsilon_{2},-2s\delta\}\cr
&
\widetilde\Pi^{'}=\{2t_{2}\varepsilon_{2},-s\delta+t_{1}\varepsilon_{1}-t_{2}\varepsilon_{2},2s\delta\},
\end{align*}
where $s,t_{1},t_{2}\in \{1,-1\}$.

\subsection{The Weyl group and Coxeter graph}
Recall that $\ggg_\bz$ admits the Weyl group $\langle \widehat r_{2\delta}\rangle\times  \langle \widehat r_{2\varepsilon_1}\rangle\times \langle \widehat r_{2\varepsilon_2}\rangle$ which is isomorphic to $\bbz_2\times \bbz_2\times \bbz_2$. Now we can calculate the Coxeter graph for  super Weyl group of $D(2,1,\alpha)$.

\begin{theorem}\label{thm: excep d}
	The Coxeter graph of $\wsc$ for $D(2,1,\alpha)$  is	

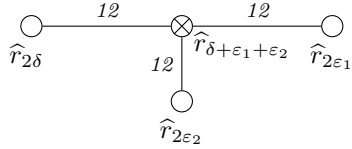
\begin{center}
	\begin{tikzpicture}
		
		\node[circle,draw,minimum size=8pt,inner sep=0pt,fill=white] (r1) at (2,0) {};

		\node[circle,draw,minimum size=8pt,inner sep=0pt,fill=white] (r2) at (0,1) {};

		\node[circle,draw,minimum size=8pt,inner sep=0pt,fill=white] (r3) at (2,1) {};
		\draw (r3.135) -- (r3.315);
		\draw(r3.225) -- (r3.45);
		
		\node[circle,draw,minimum size=8pt,inner sep=0pt,fill=white] (r4) at (4,1) {};

		\draw (r2) -- node[above,sloped] {\footnotesize 12} (r3);
		\draw (r1) -- node[left] {\footnotesize 12} (r3);
		\draw (r3) -- node[above,sloped] {\footnotesize 12} (r4);
		
		\node[below=3pt] at (r1) {$\widehat r_{2\varepsilon_2} $};
		\node[below=3pt] at (r2) {$ \widehat r_{2\delta~~~~}$~~~~~~~};
		\node at(2.8,0.75){ $\widehat r_{\delta+\varepsilon_{1}+\varepsilon_2}$};
		\node[below=3pt] at (r4) {$\widehat r_{2\varepsilon_{1}}$};
		
	\end{tikzpicture}
\end{center}

\end{theorem}

 Note that the Coxeter graph of  the subsystem $(\scrw, \spo)$ is already known, which is just the one of three vertices mutually isolated.
	We only need to calculate  the order of  $\widehat r_{\delta+\varepsilon_{1}+\varepsilon_{2}}\cdot \widehat r_{2\delta}$, $\widehat r_{\delta+\varepsilon_{1}+\varepsilon_{2}}\cdot \widehat r_{2\varepsilon_{1}}$ and $\widehat r_{\delta+\varepsilon_{1}+\varepsilon_{2}}\cdot \widehat r_{2\varepsilon_{2}}$.
	
	First, we need to  calculate the actions of  $\widehat r_{\delta+\varepsilon_{1}+\varepsilon_{2}}\cdot \widehat r_{2\delta}$ on all the  fundamental systems, and get the order of $\widehat r_{\delta+\varepsilon_{1}+\varepsilon_{2}}\cdot \widehat r_{2\delta}$ on each fundamental system.  The least common multiple of these  orders is $12$. So the order of $\widehat r_{\delta+\varepsilon_{1}+\varepsilon_{2}}\cdot \widehat r_{2\delta}$ is $12$.  Denote $\widehat r_{\delta12}$ by $\widehat r_{\delta+\varepsilon_{1}+\varepsilon_{2}}$.
{We can discuss the order of $\widehat r_{\delta12}\cdot \widehat r_{2\delta}$  on different fundamental systems in following cases:

(i) $\widetilde\Pi^{'}$ contains elements $\pm(-\delta+\varepsilon_{1}+\varepsilon_{2})$ and $\pm(\delta+\varepsilon_{1}+\varepsilon_{2})$. By observing the fundamental systems, we can see that $\widetilde\Pi^{'}=\{s(\delta+\varepsilon_{1}-\varepsilon_{2}),s(-\delta+\varepsilon_{1}+\varepsilon_{2}),s(-\delta-\varepsilon_{1}-\varepsilon_{2})\}$ or $\{s(\delta-\varepsilon_{1}+\varepsilon_{2}),s(-\delta-\varepsilon_{1}-\varepsilon_{2}),s(-\delta+\varepsilon_{1}+\varepsilon_{2})\}$, where $s\in \{\pm1\}$. By calculation, we have
\begin{align*}
&\{s(\delta+\varepsilon_{1}-\varepsilon_{2}),s(-\delta+\varepsilon_{1}+\varepsilon_{2}),s(-\delta-\varepsilon_{1}-\varepsilon_{2})\}\overset
	{\widehat r_{2\delta}}{\longrightarrow}\{s(-\delta+\varepsilon_{1}-\varepsilon_{2}),s(\delta+\varepsilon_{1}+\varepsilon_{2}),s(\delta-\varepsilon_{1}-\varepsilon_{2})\}\overset
	{\widehat r_{\delta12}}{\longrightarrow}\cr
&\{s(2\varepsilon_{1}),s(-\delta-\varepsilon_{1}-\varepsilon_{2}),s(2\delta)\} \overset
	{\widehat r_{2\delta}}{\longrightarrow}\{s(2\varepsilon_{1}),s(\delta-\varepsilon_{1}-\varepsilon_{2}),s(-2\delta)\}\overset
	{\widehat r_{\delta12}}{\longrightarrow} \cr &\{s(2\varepsilon_{1}),s(\delta-\varepsilon_{1}-\varepsilon_{2}),s(-2\delta)\}\overset
	{\widehat r_{2\delta}}{\longrightarrow} \{s(2\varepsilon_{1}),s(-\delta-\varepsilon_{1}-\varepsilon_{2}),s(2\delta)\}\overset
	{\widehat r_{\delta12}}{\longrightarrow}\cr
 &\{s(-\delta+\varepsilon_{1}-\varepsilon_{2}),s(\delta+\varepsilon_{1}+\varepsilon_{2}),s(\delta-\varepsilon_{1}-\varepsilon_{2})\}\overset
	{\widehat r_{2\delta}}{\longrightarrow} \{s(\delta+\varepsilon_{1}-\varepsilon_{2}),s(-\delta+\varepsilon_{1}+\varepsilon_{2}),s(-\delta-\varepsilon_{1}-\varepsilon_{2})\}\overset
	{\widehat r_{\delta12}}{\longrightarrow}\cr
 &\{ s(-2\varepsilon_{2}),s(-2\delta),s(\delta+\varepsilon_{1}+\varepsilon_{2})\}\overset
 {\widehat r_{2\delta}}{\longrightarrow}
 \{ s(-2\varepsilon_{2}),s(2\delta),s(-\delta+\varepsilon_{1}+\varepsilon_{2})\}\overset
 {\widehat r_{\delta12}}{\longrightarrow}\cr
 &\{ s(-2\varepsilon_{2}),s(2\delta),s(-\delta+\varepsilon_{1}+\varepsilon_{2})\}\overset
 {\widehat r_{2\delta}}{\longrightarrow}
 \{ s(-2\varepsilon_{2}),s(-2\delta),s(\delta+\varepsilon_{1}+\varepsilon_{2})\}
 \overset
 {\widehat r_{\delta12}}{\longrightarrow}\cr
 &\{s(\delta+\varepsilon_{1}-\varepsilon_{2}),s(-\delta+\varepsilon_{1}+\varepsilon_{2}),s(-\delta-\varepsilon_{1}-\varepsilon_{2})\}.
\end{align*}
So the order of the action $\widehat r_{\delta12}\cdot \widehat r_{2\delta}$ on $\widetilde\Pi^{'}=\{s(\delta+\varepsilon_{1}-\varepsilon_{2}),s(-\delta+\varepsilon_{1}+\varepsilon_{2}),s(-\delta-\varepsilon_{1}-\varepsilon_{2})\}$ is $6$. Similarly, the order of the action $\widehat r_{\delta12}\cdot \widehat r_{2\delta}$ on $\widetilde\Pi^{'}=\{s(\delta-\varepsilon_{1}+\varepsilon_{2}),s(-\delta-\varepsilon_{1}-\varepsilon_{2}),s(-\delta+\varepsilon_{1}+\varepsilon_{2})\}$ is also $6$, which is left to readers.
}	

(ii) $\widetilde\Pi^{'}$ contains elements $\pm(\delta+\varepsilon_{1}+\varepsilon_{2})$ and $\pm2\delta$, then $\widetilde\Pi^{'}=\{s(\delta+\varepsilon_{1}+\varepsilon_{2}),s(-2\delta),s(-2\varepsilon_{i})\}$ $(i=1,2)$.
\begin{align*}
	&\{s(\delta+\varepsilon_{1}+\varepsilon_{2}),s(-2\delta),s(-2\varepsilon_{i})\}\overset
	{\widehat r_{2\delta}}{\longrightarrow}\{s(-\delta+\varepsilon_{1}+\varepsilon_{2}),s(2\delta),s(-2\varepsilon_{i})\}\overset
	{\widehat r_{\delta12}}{\longrightarrow}\cr
	&\{s(-\delta+\varepsilon_{1}+\varepsilon_{2}),s(2\delta),s(-2\varepsilon_{i})\} \overset
	{\widehat r_{2\delta}}{\longrightarrow}\{s(\delta+\varepsilon_{1}+\varepsilon_{2}),s(-2\delta),s(-2\varepsilon_{i})\}\overset
	{\widehat r_{\delta12}}{\longrightarrow} \cr &\{s(-\delta-\varepsilon_{1}-\varepsilon_{2}),s(-\delta+\varepsilon_{1}+\varepsilon_{2}),s(\delta+\varepsilon_{1}+\varepsilon_{2}-2\varepsilon_{i})\}\overset
	{\widehat r_{2\delta}}{\longrightarrow} \{s(\delta-\varepsilon_{1}-\varepsilon_{2}),s(\delta+\varepsilon_{1}+\varepsilon_{2}),s(-\delta+\varepsilon_{1}+\varepsilon_{2}-2\varepsilon_{i})\}\overset
	{\widehat r_{\delta12}}{\longrightarrow}\cr
	&\{s(2\delta),s(-\delta-\varepsilon_{1}-\varepsilon_{2}),s(2\varepsilon_{1}+2\varepsilon_{2}-2\varepsilon_{i})\}\overset
	{\widehat r_{2\delta}}{\longrightarrow} \{s(-2\delta),s(\delta-\varepsilon_{1}-\varepsilon_{2}),s(2\varepsilon_{1}+2\varepsilon_{2}-2\varepsilon_{i})\}\overset
	{\widehat r_{\delta12}}{\longrightarrow}\cr
	&\{s(-2\delta),s(\delta-\varepsilon_{1}-\varepsilon_{2}),s(2\varepsilon_{1}+2\varepsilon_{2}-2\varepsilon_{i})\}\overset
	{\widehat r_{2\delta}}{\longrightarrow}
	\{s(2\delta),s(-\delta-\varepsilon_{1}-\varepsilon_{2}),s(2\varepsilon_{1}+2\varepsilon_{2}-2\varepsilon_{i})\}\overset
	{\widehat r_{\delta12}}{\longrightarrow}\cr
	&\{s(\delta-\varepsilon_{1}-\varepsilon_{2}),s(\delta+\varepsilon_{1}+\varepsilon_{2}),s(-\delta+\varepsilon_{1}+\varepsilon_{2}-2\varepsilon_{i})\}\overset
	{\widehat r_{2\delta}}{\longrightarrow}
	\{s(-\delta-\varepsilon_{1}-\varepsilon_{2}),s(-\delta+\varepsilon_{1}+\varepsilon_{2}),s(\delta+\varepsilon_{1}+\varepsilon_{2}-2\varepsilon_{i})\}
	\overset
	{\widehat r_{\delta12}}{\longrightarrow}\cr
	&\{s(\delta+\varepsilon_{1}+\varepsilon_{2}),s(-2\delta),s(-2\varepsilon_{i})\}.
\end{align*}
So the order of the action of  $\widehat r_{\delta+\varepsilon_{1}+\varepsilon_{2}}\cdot \widehat r_{2\delta}$ on this fundamental system is $6$.

(iii) $\widetilde\Pi^{'}$ contains elements $\pm(-\delta+\varepsilon_{1}+\varepsilon_{2})$ and $\pm(2\delta)$, then $\widetilde\Pi^{'}=\{s(-\delta+\varepsilon_{1}+\varepsilon_{2}),s(2\delta),s(-2\varepsilon_{i})\}$ $(i=1,2)$. Similar to the case  $\widetilde\Pi^{'}=\{s(\delta+\varepsilon_{1}+\varepsilon_{2}),s(-2\delta),s(-2\varepsilon_{i})\}$ $(i=1,2)$, one can easily check that the order in this case is also $6$.

(iv) $\widetilde\Pi^{'}=\{s(\delta+\varepsilon_{1}+\varepsilon_{2}),\alpha,\beta\}$, where $\alpha,\beta\ne \pm(\delta-\varepsilon_{1}-\varepsilon_{2}),\pm(2\delta)$. By observing all fundamental systems, one can easily see that $\widehat r_{2\delta}(\alpha),\widehat r_{2\delta}(\beta),\widehat r_{2\delta}\cdot \widehat r_{\delta12}(\alpha),\widehat r_{2\delta}\cdot \widehat r_{\delta12}(\beta)\ne \pm(\delta+\varepsilon_{1}+\varepsilon_{2})$. Then it follows that
\begin{align*}
	&\{s(\delta+\varepsilon_{1}+\varepsilon_{2}),\alpha,\beta\}\overset
	{\widehat r_{2\delta}}{\longrightarrow}\{s(-\delta+\varepsilon_{1}+\varepsilon_{2}),\widehat r_{2\delta}(\alpha),\widehat r_{2\delta}(\beta)\}\overset
	{\widehat r_{\delta12}}{\longrightarrow}\cr
	&\{s(-\delta+\varepsilon_{1}+\varepsilon_{2}),\widehat r_{2\delta}(\alpha),\widehat r_{2\delta}(\beta)\} \overset
	{\widehat r_{2\delta}}{\longrightarrow}\{s(\delta+\varepsilon_{1}+\varepsilon_{2}),\alpha,\beta\}\overset
	{\widehat r_{\delta12}}{\longrightarrow} \cr &\{s(-\delta-\varepsilon_{1}-\varepsilon_{2}),\widehat r_{\delta12}(\alpha),\widehat r_{\delta12}(\beta)\}\overset
	{\widehat r_{2\delta}}{\longrightarrow} \{s(\delta-\varepsilon_{1}-\varepsilon_{2}),\widehat r_{2\delta}\cdot \widehat r_{\delta12}(\alpha),\widehat r_{2\delta}\cdot \widehat r_{\delta12}(\beta)\}\overset
	{\widehat r_{\delta12}}{\longrightarrow}\cr
	&\{s(\delta-\varepsilon_{1}-\varepsilon_{2}),\widehat r_{2\delta}\cdot \widehat r_{\delta12}(\alpha),\widehat r_{2\delta}\cdot \widehat r_{\delta12}(\beta)\}\overset
	{\widehat r_{2\delta}}{\longrightarrow} \{s(-\delta-\varepsilon_{1}-\varepsilon_{2}),\widehat r_{\delta12}(\alpha),\widehat r_{\delta12}(\beta)\}\overset
	{\widehat r_{\delta12}}{\longrightarrow}\cr
	&\{s(\delta+\varepsilon_{1}+\varepsilon_{2}),\alpha,\beta\}.
\end{align*}
Hence, the order of $\widehat r_{\delta12}\cdot \widehat r_{2\delta}$ in this case is $4$. Similarly, one can check the following case (v) and we leave the proof to readers.

(v) $\widetilde\Pi^{'}=\{s(-\delta+\varepsilon_{1}+\varepsilon_{2}),\alpha,\beta\}$, where $\alpha,\beta\ne \pm(\delta+\varepsilon_{1}+\varepsilon_{2}),\pm(2\delta)$. The order of $\widehat r_{\delta12}\cdot \widehat r_{2\delta}$ on these fundamental systems is $4$.

(vi) $\widetilde\Pi^{'}$ contains no elements $\pm(\pm\delta+\varepsilon_{1}+\varepsilon_{2})$, then in this case all elements in fundamental systems obtained by repeatedly acting real reflection $\widehat r_{2\delta}$ and odd reflection $\widehat r_{\delta12}$ on $\widetilde\Pi^{'}$ are fixed by $\widehat r_{\delta12}$. Thus the order is 2.

Summing up, we conclude that the order of $\widehat r_{\delta+\varepsilon_{1}+\varepsilon_{2}}\cdot \widehat r_{2\delta}$ is $12$.

Next we will calculate the order of $\widehat r_{\delta12}\cdot \widehat r_{2\varepsilon_{1}}$. We also discuss the order of $\widehat r_{\delta12}\cdot \widehat r_{2\varepsilon_{1}}$ on different fundamental systems as follow:

(i) $\widetilde\Pi^{'}$ contains elements $\pm(\delta\pm\varepsilon_{1}+\varepsilon_{2})$ and $\pm2\varepsilon_{1}$, then $\widetilde\Pi^{'}=\{s(\delta+\varepsilon_{1}+\varepsilon_{2}),s(-2\varepsilon_{1}),\alpha\}$ or $\{s(\delta-\varepsilon_{1}+\varepsilon_{2}),s(2\varepsilon_{1}),\alpha\}$. Here we just give the calculation of $\widetilde\Pi^{'}=\{s(\delta+\varepsilon_{1}+\varepsilon_{2}),s(-2\varepsilon_{1}),\alpha\}$ and the case  $\widetilde\Pi^{'}=\{s(\delta-\varepsilon_{1}+\varepsilon_{2}),s(2\varepsilon_{1}),\alpha\}$ is similar.
For $\widetilde\Pi^{'}=\{s(\delta+\varepsilon_{1}+\varepsilon_{2}),s(-2\varepsilon_{1}),\alpha\}$, one can see that $\alpha=-2\varepsilon_{2}$ or $-2\delta$. For both cases, we have $\widehat r_{2\varepsilon_{1}}(\alpha),(\widehat r_{2\varepsilon_{1}}\cdot\widehat r_{\delta12})^{2}(\alpha)\ne\pm(\delta+\varepsilon_{1}+\varepsilon_{2})$. Then we get
\begin{align*}
	&\{s(\delta+\varepsilon_{1}+\varepsilon_{2}),s(-2\varepsilon_{1}),\alpha\}\overset
	{\widehat r_{2\varepsilon_{1}}}{\longrightarrow}\{s(\delta-\varepsilon_{1}+\varepsilon_{2}),s(2\varepsilon_{1}),\widehat r_{2\varepsilon_{1}}(\alpha)\}\overset
	{\widehat r_{\delta12}}{\longrightarrow}\cr
	&\{s(\delta-\varepsilon_{1}+\varepsilon_{2}),s(2\varepsilon_{1}),\widehat r_{2\varepsilon_{1}}(\alpha)\} \overset
	{\widehat r_{2\varepsilon_{1}}}{\longrightarrow}\{s(\delta+\varepsilon_{1}+\varepsilon_{2}),s(-2\varepsilon_{1}),\alpha\}\overset
	{\widehat r_{\delta12}}{\longrightarrow} \cr &\{s(-\delta-\varepsilon_{1}-\varepsilon_{2}),s(\delta-\varepsilon_{1}+\varepsilon_{2}),\widehat r_{\delta12}(\alpha)\}\overset
	{\widehat r_{2\varepsilon_{1}}}{\longrightarrow} \{s(-\delta+\varepsilon_{1}-\varepsilon_{2}),s(\delta+\varepsilon_{1}+\varepsilon_{2}),\widehat r_{2\varepsilon_{1}}\cdot \widehat r_{\delta12}(\alpha)\}\overset
	{\widehat r_{\delta12}}{\longrightarrow}\cr
	&\{s(2\varepsilon_{1}),s(-\delta-\varepsilon_{1}-\varepsilon_{2}),\widehat r_{\delta12}\cdot\widehat r_{2\varepsilon_{1}}\cdot \widehat r_{\delta12}(\alpha)\}\overset
	{\widehat r_{2\varepsilon_{1}}}{\longrightarrow} \{s(-2\varepsilon_{1}),s(-\delta+\varepsilon_{1}-\varepsilon_{2}),(\widehat r_{2\varepsilon_{1}}\cdot\widehat r_{\delta12})^{2}(\alpha)\}\overset
	{\widehat r_{\delta12}}{\longrightarrow}\cr
	&\{s(-2\varepsilon_{1}),s(-\delta+\varepsilon_{1}-\varepsilon_{2}),(\widehat r_{2\varepsilon_{1}}\cdot\widehat r_{\delta12})^{2}(\alpha)\}\overset
	{\widehat r_{2\varepsilon_{1}}}{\longrightarrow}\{s(2\varepsilon_{1}),s(-\delta-\varepsilon_{1}-\varepsilon_{2}),\widehat r_{\delta12}\cdot\widehat r_{2\varepsilon_{1}}\cdot \widehat r_{\delta12}(\alpha)\}\overset
	{\widehat r_{\delta12}}{\longrightarrow}\cr
	&\{s(-\delta+\varepsilon_{1}-\varepsilon_{2}),s(\delta+\varepsilon_{1}+\varepsilon_{2}),\widehat r_{2\varepsilon_{1}}\cdot \widehat r_{\delta12}(\alpha)\}\overset
	{\widehat r_{2\varepsilon_{1}}}{\longrightarrow}\{s(-\delta-\varepsilon_{1}-\varepsilon_{2}),s(\delta-\varepsilon_{1}+\varepsilon_{2}),\widehat r_{\delta12}(\alpha)\}\}\overset
	{\widehat r_{\delta12}}{\longrightarrow}\cr
	&\{s(\delta+\varepsilon_{1}+\varepsilon_{2}),s(-2\varepsilon_{1}),\alpha\}
\end{align*}
So in this case the order is $6$.

(ii) $\widetilde\Pi^{'}$ contains elements $\pm(\delta+\varepsilon_{1}+\varepsilon_{2})$ and $\pm(\delta-\varepsilon_{1}+\varepsilon_{2})$, then $\widetilde\Pi^{'}=\{s(\delta+\varepsilon_{1}+\varepsilon_{2}),s(-\delta+\varepsilon_{1}-\varepsilon_{2}),\alpha\}$, where $\alpha=s(-\delta-\varepsilon_{1}+\varepsilon_{2})$ or $s(\delta-\varepsilon_{1}-\varepsilon_{2})$. For both cases, $\widehat r_{2\varepsilon_{1}}\cdot\widehat r_{\delta12}\cdot \widehat r_{2\varepsilon_{1}}(\alpha),\widehat r_{2\varepsilon_{1}}\cdot\widehat r_{\delta12}(\alpha)\ne\pm(\delta+\varepsilon_{1}+\varepsilon_{2})$. Then we have
 \begin{align*}
 	&\{s(\delta+\varepsilon_{1}+\varepsilon_{2}),s(-\delta+\varepsilon_{1}-\varepsilon_{2}),\alpha\}\overset
 	{\widehat r_{2\varepsilon_{1}}}{\longrightarrow}\{s(\delta-\varepsilon_{1}+\varepsilon_{2}),s(-\delta-\varepsilon_{1}-\varepsilon_{2}),\widehat r_{2\varepsilon_{1}}(\alpha)\}\overset
 	{\widehat r_{\delta12}}{\longrightarrow}\cr
 	&\{s(-2\varepsilon_{1}),s(\delta+\varepsilon_{1}+\varepsilon_{2}),\widehat r_{\delta12}\cdot \widehat r_{2\varepsilon_{1}}(\alpha)\} \overset
 	{\widehat r_{2\varepsilon_{1}}}{\longrightarrow}\{s(2\varepsilon_{1}),s(\delta-\varepsilon_{1}+\varepsilon_{2}),\widehat r_{2\varepsilon_{1}}\cdot\widehat r_{\delta12}\cdot \widehat r_{2\varepsilon_{1}}(\alpha)\}\overset
 	{\widehat r_{\delta12}}{\longrightarrow} \cr &\{s(2\varepsilon_{1}),s(\delta-\varepsilon_{1}+\varepsilon_{2}),\widehat r_{2\varepsilon_{1}}\cdot\widehat r_{\delta12}\cdot \widehat r_{2\varepsilon_{1}}(\alpha)\}\overset
 	{\widehat r_{2\varepsilon_{1}}}{\longrightarrow} \{s(-2\varepsilon_{1}),s(\delta+\varepsilon_{1}+\varepsilon_{2}),\widehat r_{\delta12}\cdot \widehat r_{2\varepsilon_{1}}(\alpha)\}\overset
 	{\widehat r_{\delta12}}{\longrightarrow}\cr
 	&\{s(\delta-\varepsilon_{1}+\varepsilon_{2}),s(-\delta-\varepsilon_{1}-\varepsilon_{2}),\widehat r_{2\varepsilon_{1}}(\alpha)\}\overset
 	{\widehat r_{2\varepsilon_{1}}}{\longrightarrow} \{s(\delta+\varepsilon_{1}+\varepsilon_{2}),s(-\delta+\varepsilon_{1}-\varepsilon_{2}),\alpha\}\overset
 	{\widehat r_{\delta12}}{\longrightarrow}\cr
 	&\{s(-\delta-\varepsilon_{1}-\varepsilon_{2}),s(2\varepsilon_{1}),\widehat r_{\delta12}(\alpha)\}\overset
 	{\widehat r_{2\varepsilon_{1}}}{\longrightarrow}\{s(-\delta+\varepsilon_{1}-\varepsilon_{2}),s(-2\varepsilon_{1}),\widehat r_{2\varepsilon_{1}}\cdot\widehat r_{\delta12}(\alpha)\}\overset
 	{\widehat r_{\delta12}}{\longrightarrow}\cr
 	&\{s(-\delta+\varepsilon_{1}-\varepsilon_{2}),s(-2\varepsilon_{1}),\widehat r_{2\varepsilon_{1}}\cdot\widehat r_{\delta12}(\alpha)\}\overset
 	{\widehat r_{2\varepsilon_{1}}}{\longrightarrow}\{s(-\delta-\varepsilon_{1}-\varepsilon_{2}),s(2\varepsilon_{1}),\widehat r_{\delta12}(\alpha)\}\overset
 	{\widehat r_{\delta12}}{\longrightarrow}\cr
 	&\{s(\delta+\varepsilon_{1}+\varepsilon_{2}),s(-\delta+\varepsilon_{1}-\varepsilon_{2}),\alpha\}.
 \end{align*}
In this case, the order is $6$.

(iii) $\widetilde\Pi^{'}=\{s(\delta+\varepsilon_{1}+\varepsilon_{2}),\alpha,\beta\}$, where $\alpha,\beta\ne\pm(2\varepsilon_{1}),\pm(\delta-\varepsilon_{1}+\varepsilon_{2})$. By observation, we have $\widehat r_{2\varepsilon_{1}}(\alpha),\widehat r_{2\varepsilon_{1}}(\beta),\widehat r_{2\varepsilon_{1}}\cdot\widehat r_{\delta12}(\alpha),\widehat r_{2\varepsilon_{1}}\cdot\widehat r_{\delta12}(\beta)\ne\pm(\delta+\varepsilon_{1}+\varepsilon_{2})$. By calculation, we have
\begin{align*}
	&\{s(\delta+\varepsilon_{1}+\varepsilon_{2}),\alpha,\beta\}\overset
	{\widehat r_{2\varepsilon_{1}}}{\longrightarrow}\{s(\delta-\varepsilon_{1}+\varepsilon_{2}),\widehat r_{2\varepsilon_{1}}(\alpha),\widehat r_{2\varepsilon_{1}}(\beta)\}\overset
	{\widehat r_{\delta12}}{\longrightarrow}\cr
	&\{s(\delta-\varepsilon_{1}+\varepsilon_{2}),\widehat r_{2\varepsilon_{1}}(\alpha),\widehat r_{2\varepsilon_{1}}(\beta)\} \overset
	{\widehat r_{2\varepsilon_{1}}}{\longrightarrow}\{s(\delta+\varepsilon_{1}+\varepsilon_{2}),\alpha,\beta\}\overset
	{\widehat r_{\delta12}}{\longrightarrow} \cr &\{s(-\delta-\varepsilon_{1}-\varepsilon_{2}),\widehat r_{\delta12}(\alpha),\widehat r_{\delta12}(\beta)\}\overset
	{\widehat r_{2\varepsilon_{1}}}{\longrightarrow} \{s(-\delta+\varepsilon_{1}-\varepsilon_{2}),\widehat r_{2\varepsilon_{1}}\cdot\widehat r_{\delta12} (\alpha),\widehat r_{2\varepsilon_{1}}\cdot\widehat r_{\delta12} (\beta)\}\overset
	{\widehat r_{\delta12}}{\longrightarrow}\cr
	&\{s(-\delta+\varepsilon_{1}-\varepsilon_{2}),\widehat r_{2\varepsilon_{1}}\cdot\widehat r_{\delta12} (\alpha),\widehat r_{2\varepsilon_{1}}\cdot\widehat r_{\delta12} (\beta)\}\overset
	{\widehat r_{2\varepsilon_{1}}}{\longrightarrow} \{s(-\delta-\varepsilon_{1}-\varepsilon_{2}),\widehat r_{\delta12} (\alpha),\widehat r_{\delta12} (\beta)\}\overset
	{\widehat r_{\delta12}}{\longrightarrow}\cr
	&\{s(\delta+\varepsilon_{1}+\varepsilon_{2}),\alpha,\beta\}.
\end{align*}
In this case, the order is $4$.

Similarly, we can get the order of case (iv), which is also $4$ and we omit the details.

(iv) $\widetilde\Pi^{'}=\{s(\delta-\varepsilon_{1}+\varepsilon_{2}),\alpha,\beta\}$, where $\alpha,\beta\ne\pm(2\varepsilon_{1}),\pm(\delta+\varepsilon_{1}+\varepsilon_{2})$.

(v) $\widetilde\Pi^{'}$ contains no elements $\pm(\delta+\varepsilon_{1}+\varepsilon_{2})$ and $\pm(\delta-\varepsilon_{1}+\varepsilon_{2})$. The action of the odd reflection is trivial and the order in this case is $2$.

Finally, we conclude that the order of $\widehat r_{\delta12}\cdot \widehat r_{2\varepsilon_{1}}$ is $12$.

Since the weights $\varepsilon_{1}$ and $\varepsilon_{2}$ play the same role for $D(2,1,\alpha)$, we can get the order of $\widehat r_{\delta12}\cdot \widehat r_{2\varepsilon_{2}}$ is also $12$.

Now we finish the calculations.


\section {The root system and Coxeter graph for $F(4)$}
\subsection{Root system for $F(4)$}{
Let $\ggg=F(4)$. Recall that
$\ggg_\bz\cong \mathfrak{sl}(2)\oplus \mathfrak{so}(7)$ with Cartan subalgebra $\hhh$ and $\ggg_\bo=\bbc^2\oplus \bbc^8$, as a $\ggg_\bz$-module. Here $\bbc^8$ is  the $8$-dimensional spin representation of $\mathfrak{so}(7)$. The root system of $\ggg$ can be
described via the basis  $\{\varepsilon_{1},\varepsilon_{2}, \varepsilon_{3}, \delta\}$ in $\hhh^*$ with $(\delta,\delta)=-3, (\varepsilon_{i},\varepsilon_{i})=1, (\varepsilon_{i},\varepsilon_{j})=0, (\delta,\varepsilon_{i})=0, i,j =1,2,3, i\neq j$.
Moreover,
$\Phi=\Phi_{\bar{0}}\bigcup \Phi_{\bar{1}}$ is precisely presented as
$$\{\pm\varepsilon_{i}\pm\varepsilon_{j},\pm\varepsilon_{i},\pm\delta\}\cup
\{\frac{1}{2}(\pm\varepsilon_{1}\pm\varepsilon_{2}
\pm\varepsilon_{3}\pm\delta)\}.$$


      \subsection{Coxeter graph} We can take $$\Pi=\{\frac{1}{2}(\delta+\varepsilon_1+ \varepsilon_2+\varepsilon_3), -\varepsilon_3, \varepsilon_3-\varepsilon_1,  \varepsilon_1-\varepsilon_2\}$$
 as  the standard fundamental system of $F(4)$. Then we can take
  $$\widetilde \Pi=\{-\delta\}\cup \{\frac{1}{2}(\delta+\varepsilon_1+ \varepsilon_2+\varepsilon_3), -\varepsilon_3, \varepsilon_3-\varepsilon_1,  \varepsilon_1-\varepsilon_2\}$$}
as the extended standard fundamental system.

All fundamental systems can be classified into following $6$ cases:
\begin{align*}
	& \widetilde\Pi^{'}=\{\frac{1}{2}(s\delta+t_{i}\varepsilon_{i}+t_{j}\varepsilon_{j}+t_{k}\varepsilon_{k}),-t_{i}\varepsilon_{i},t_{i}\varepsilon_{i}-t_{j}\varepsilon_{j},t_{j}\varepsilon_{j}-t_{k}\varepsilon_{k}\} \cr
	&\widetilde\Pi^{'}=\{\frac{1}{2}(s\delta-t_{i}\varepsilon_{i}-t_{j}\varepsilon_{j}+t_{k}\varepsilon_{k}),\frac{1}{2}(-s\delta+t_{i}\varepsilon_{i}+t_{j}\varepsilon_{j}+t_{k}\varepsilon_{k}),t_{i}\varepsilon_{i}-t_{j}\varepsilon_{j},t_{j}\varepsilon_{j}-t_{k}\varepsilon_{k}\} \cr
	&\widetilde\Pi^{'}=\{\frac{1}{2}(s\delta-t_{i}\varepsilon_{i}+t_{j}\varepsilon_{j}-t_{k}\varepsilon_{k}),\frac{1}{2}(-s\delta+t_{i}\varepsilon_{i}+t_{j}\varepsilon_{j}-t_{k}\varepsilon_{k}),t_{i}\varepsilon_{i}-t_{j}\varepsilon_{j},t_{k}\varepsilon_{k}\} \cr
	&\widetilde\Pi^{'}=\{\frac{1}{2}(s\delta+t_{i}\varepsilon_{i}-t_{j}\varepsilon_{j}-t_{k}\varepsilon_{k}),\frac{1}{2}(-s\delta+t_{i}\varepsilon_{i}-t_{j}\varepsilon_{j}+t_{k}\varepsilon_{k}),t_{i}\varepsilon_{i}-t_{j}\varepsilon_{j},\frac{1}{2}(s\delta-t_{i}\varepsilon_{i}+t_{j}\varepsilon_{j}+t_{k}\varepsilon_{k})\} \cr
	&\widetilde\Pi^{'}=\{\frac{1}{2}(-s\delta+t_{i}\varepsilon_{i}-t_{j}\varepsilon_{j}-t_{k}\varepsilon_{k}),s\delta,t_{j}\varepsilon_{j}-t_{k}\varepsilon_{k},t_{k}\varepsilon_{k}\} \cr
	&\widetilde\Pi^{'}=\{\frac{1}{2}(-s\delta-t_{i}\varepsilon_{i}+t_{j}\varepsilon_{j}+t_{k}\varepsilon_{k}),s\delta,t_{i}\varepsilon_{i}-t_{j}\varepsilon_{j},t_{j}\varepsilon_{j}-t_{k}\varepsilon_{k}\},
\end{align*}
where $t_{i},t_{j},t_{k}\in \{1,-1\}$, $\{i,j,k\}=\{1,2,3\}$.


\begin{theorem}\label{thm: excep f}
	The Coxeter graph of $\wsc$ for $F(4)$ is	

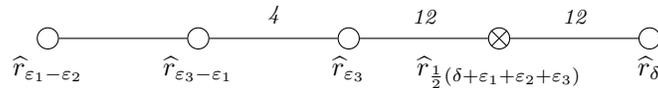
\begin{center}
	\begin{tikzpicture}
		
		\node[circle,draw,minimum size=8pt,inner sep=0pt,fill=white] (r1) at (-5,1) {};

		\node[circle,draw,minimum size=8pt,inner sep=0pt,fill=white] (r2) at (-3,1) {};
		
		\node[circle,draw,minimum size=8pt,inner sep=0pt,fill=white] (r3) at (-1,1) {};
		
		\node[circle,draw,minimum size=8pt,inner sep=0pt,fill=white] (r4) at (1,1) {};
			\draw (r4.135) -- (r4.315);
		\draw(r4.225) -- (r4.45);
		
		\node[circle,draw,minimum size=8pt,inner sep=0pt,fill=white] (r5) at (3,1) {};
		
		\draw (r1) -- (r2);
		\draw (r2) -- node[above=2pt] {\footnotesize 4} (r3);
		\draw (r3) -- node[above=2pt] {\footnotesize 12} (r4);
		\draw (r4) -- node[above=2pt] {\footnotesize 12} (r5);
		
		\node[below=3pt] at (r1) {$\widehat r_{\varepsilon_1-\varepsilon_2} $};
		\node[below=3pt] at (r2) {$ \widehat r_{\varepsilon_3-\varepsilon_1}$};
		\node[below=3pt] at (r3) {$\widehat r_{\varepsilon_3}$};
		\node[below=3pt] at (r4) {$\widehat r_{\frac{1}{2}(\delta+\varepsilon_{1}+\varepsilon_2+\varepsilon_3)}$};		
		\node[below=3pt] at (r5) {$\widehat r_{\delta}$};
		
	\end{tikzpicture}
\end{center}

\end{theorem}

 {
Note that the Coxeter graph of  the subsystem $(\scrw, \spo)$ is already known (see for example, \cite[Theorem 1.3.3]{GP}). It suffices for us to compute $m_{xy}$ involving $\widehat r_{\frac{1}{2}(\delta+\varepsilon_{1}+\varepsilon_2+\varepsilon_3)}$
.} 	
	Here we just give the calculation of $\widehat r_{\delta}\cdot\widehat r_{\frac{1}{2}(\delta+\varepsilon_{1}+\varepsilon_2+\varepsilon_3)}$. The calculation of  $\widehat r_{\varepsilon_{3}}\cdot\widehat r_{\frac{1}{2}(\delta+\varepsilon_{1}+\varepsilon_2+\varepsilon_3)}$ is  similar, if one can notice that $\widehat r_{\delta}$ and $\widehat r_{\varepsilon_{3}}$ play the same role in the action by exchanging the signs before $\delta$ and $\varepsilon_{3}$. The order of both cases is $12$.

For convenience, we denote $\gamma=\frac{1}{2}(\delta+\varepsilon_{1}+\varepsilon_2+\varepsilon_3)$. Now we need to calculate the order of $\widehat r_{\delta}\cdot \widehat r_{\gamma}$ on all fundamental systems. We also consider in different cases.

(i) $\widetilde\Pi^{'}$ contains $\pm(\gamma-\delta)$ and $\pm\delta$, then $\widetilde\Pi^{'}=\{s(\gamma-\delta),s\delta,\alpha_{1},\alpha_{2}\}$, $s\in\{\pm1\}$, where $\alpha_{i}$ satisfy the condition that $(\widehat r_{\gamma}\cdot \widehat r_{\delta})^{2}(\alpha_{i}),\widehat r_{\delta}\cdot (\widehat r_{\gamma}\cdot \widehat r_{\delta})^{2}(\alpha_{i})\ne \pm\gamma$, for $i=1,2$. Then we have
\begin{align*}
	&\{s(\gamma-\delta),s\delta,\alpha_{1},\alpha_{2}\}
	\overset{\widehat r_{\gamma}}{\longrightarrow}
	\{s(\gamma-\delta),s\delta,\alpha_{1},\alpha_{2}\}
	\overset{\widehat r_{\delta}}{\longrightarrow}\cr
	&	\{s\gamma, s(-\delta),\widehat r_{\delta}(\alpha_{1}),\widehat r_{\delta}(\alpha_{2}) \}
	\overset{\widehat r_{\gamma}}{\longrightarrow}
	\{s(-\gamma), s(\gamma-\delta),\widehat r_{\gamma}\cdot\widehat r_{\delta}(\alpha_{1}),\widehat r_{\gamma}\cdot\widehat r_{\delta}(\alpha_{2}) \}\overset{\widehat r_{\delta}}{\longrightarrow}\cr
	&\{s(-\gamma+\delta), s\gamma,\widehat r_{\delta}\cdot \widehat r_{\gamma}\cdot\widehat r_{\delta}(\alpha_{1}),\widehat r_{\delta}\cdot\widehat r_{\gamma}\cdot\widehat r_{\delta}(\alpha_{2}) \}
	\overset{\widehat r_{\gamma}}{\longrightarrow}
	\{s\delta, s(-\gamma),(\widehat r_{\gamma}\cdot \widehat r_{\delta})^{2}(\alpha_{1}),(\widehat r_{\gamma}\cdot \widehat r_{\delta})^{2}(\alpha_{2})\}
	\overset{\widehat r_{\delta}}{\longrightarrow}\cr
	&\{s(-\delta), s(-\gamma+\delta),\widehat r_{\delta}\cdot (\widehat r_{\gamma}\cdot \widehat r_{\delta})^{2}(\alpha_{1}),\widehat r_{\delta}\cdot (\widehat r_{\gamma}\cdot \widehat r_{\delta})^{2}(\alpha_{2})\}
	\overset{\widehat r_{\gamma}}{\longrightarrow}
	\{s(-\delta), s(-\gamma+\delta),\widehat r_{\delta}\cdot (\widehat r_{\gamma}\cdot \widehat r_{\delta})^{2}(\alpha_{1}),\widehat r_{\delta}\cdot (\widehat r_{\gamma}\cdot \widehat r_{\delta})^{2}(\alpha_{2})\}
	\overset{\widehat r_{\delta}}{\longrightarrow}\cr
	&	\{s\delta, s(-\gamma),(\widehat r_{\gamma}\cdot \widehat r_{\delta})^{2}(\alpha_{1}),(\widehat r_{\gamma}\cdot \widehat r_{\delta})^{2}(\alpha_{2}) \}\overset{\widehat r_{\gamma}}{\longrightarrow}\{s(-\gamma+\delta), s\gamma,\widehat r_{\delta}\cdot \widehat r_{\gamma}\cdot\widehat r_{\delta}(\alpha_{1}),\widehat r_{\delta}\cdot\widehat r_{\gamma}\cdot\widehat r_{\delta}(\alpha_{2}) \}\overset{\widehat r_{\delta}}{\longrightarrow}\cr
	&	\{s(-\gamma), s(\gamma-\delta),\widehat r_{\gamma}\cdot\widehat r_{\delta}(\alpha_{1}),\widehat r_{\gamma}\cdot\widehat r_{\delta}(\alpha_{2}) \}\overset{\widehat r_{\gamma}}{\longrightarrow}\{s\gamma, s(-\delta),\widehat r_{\delta}(\alpha_{1}),\widehat r_{\delta}(\alpha_{2}) \}\overset{\widehat r_{\delta}}{\longrightarrow}\cr
	&\{s(\gamma-\delta),s\delta,\alpha_{1},\alpha_{2}\}.
\end{align*}
Thus the order in this case is $6$.

(ii) $\widetilde\Pi^{'}$ contains $\pm\gamma$ and $\pm\delta$, then $\widetilde\Pi^{'}=\{s\gamma,s(-\delta),\alpha_{1},\alpha_{2}\}$, $s\in\{\pm1\}$. The order in this case is also $6$ and we leave the calculation to readers.

(iii) $\widetilde\Pi^{'}=\{s\gamma, \alpha_{1},\alpha_{2},\alpha_{3}\}$, where $\alpha_{i}\ne\pm\delta$, for $i=1,2,3$. By observation, one can easily see that $\widehat r_{\delta}(\alpha_{i}),\widehat r_{\delta}\cdot\widehat r_{\gamma}(\alpha_{i})\ne\pm\gamma$.
\begin{align*}
	&\{s\gamma,\alpha_{1},\alpha_{2},\alpha_{3}\}
	\overset{\widehat r_{\gamma}}{\longrightarrow}
	\{s(-\gamma),\widehat r_{\gamma}(\alpha_{1}),\widehat r_{\gamma}(\alpha_{2}),\widehat r_{\gamma}(\alpha_{3})\}
	\overset{\widehat r_{\delta}}{\longrightarrow}\cr
	&	\{s(-\gamma+\delta), \widehat r_{\delta}\cdot\widehat r_{\gamma}(\alpha_{1}),\widehat r_{\delta}\cdot\widehat r_{\gamma}(\alpha_{2}),\widehat r_{\delta}\cdot\widehat r_{\gamma}(\alpha_{3}) \}
	\overset{\widehat r_{\gamma}}{\longrightarrow}
	\{s(-\gamma+\delta), \widehat r_{\delta}\cdot\widehat r_{\gamma}(\alpha_{1}),\widehat r_{\delta}\cdot\widehat r_{\gamma}(\alpha_{2}),\widehat r_{\delta}\cdot\widehat r_{\gamma}(\alpha_{3}) \}\overset{\widehat r_{\delta}}{\longrightarrow}\cr
	&\{s(-\gamma),\widehat r_{\gamma}(\alpha_{1}),\widehat r_{\gamma}(\alpha_{2}),\widehat r_{\gamma}(\alpha_{3})\}
	\overset{\widehat r_{\gamma}}{\longrightarrow}
	\{s\gamma,\alpha_{1},\alpha_{2},\alpha_{3}\}
	\overset{\widehat r_{\delta}}{\longrightarrow}\cr
	&\{s(\gamma-\delta),\widehat r_{\delta}(\alpha_{1}),\widehat r_{\delta}(\alpha_{2}),\widehat r_{\delta}(\alpha_{3}) \}
	\overset{\widehat r_{\gamma}}{\longrightarrow}
	\{s(\gamma-\delta),\widehat r_{\delta}(\alpha_{1}),\widehat r_{\delta}(\alpha_{2}),\widehat r_{\delta}(\alpha_{3}) \}
	\overset{\widehat r_{\delta}}{\longrightarrow}\cr
	&	\{s\gamma,\alpha_{1},\alpha_{2},\alpha_{3}\}.
\end{align*}
Thus in this case the order is $4$.

Similarly, one can check the following case and the details are omitted.

(iv) $\widetilde\Pi^{'}=\{s(\gamma-\delta), \alpha_{1},\alpha_{2},\alpha_{3}\}$, where $\alpha_{i}\ne\pm\delta$, for $i=1,2,3$. It is easy to check the order in this case is also $4$.

(v) $\widetilde\Pi^{'}$ contains no elements $\pm\gamma,\pm\delta,\pm(\gamma-\delta)$. In this case, the action of odd reflection is trivial and the order is $2$.

Finally, we conclude that the order of $\widehat r_{\delta}\cdot\widehat r_{\frac{1}{2}(\delta+\varepsilon_{1}+\varepsilon_2+\varepsilon_3)}$ is $12$.

Then we finish the calculations.

\section {The root system and Coxeter graph for $G(3)$}
\subsection{Root system of $G(3)$} {
Let $\ggg=\ggg_\bz\oplus \ggg_\bo$ be the exceptional simple Lie superalgebra $G(3)$. We have $\ggg_\bz\cong G_2\oplus\mathfrak{sl}(2)$, and $\ggg_\bo\cong \bbc^7\oplus \bbc^2$ as an adjoint $\ggg_\bz$-module. Here $\bbc^7$ denotes the $7$-dimensional irreducible $G_2$-module, and $\bbc^2$ the natural $\mathfrak{sl}(2)$-module (see \cite{K} or \cite{CSW} for more details on $G(3)$).

The root system of $\ggg$ can be described as follows. Fix a Cartan subalgebra  $\hhh\in\ggg_\bz$. Assume $\varepsilon_i\in\hhh^*$ ($i=1,2,3$)
satisfy the linear relation $\varepsilon+\varepsilon_2+\varepsilon_3=0$.
Choose the standard simple system  $\Pi=\{\delta-\varepsilon_1,   \varepsilon_2-\varepsilon_3, -\varepsilon_2, \}$. Then the standard positive roots are $\Phi^+=\Phi^+_0\cup\Phi^+_1$
with
\begin{align}
&\Phi_0^+=\{2\delta, \varepsilon_1, -\varepsilon_2,-\varepsilon_3, \varepsilon_1-\varepsilon_2,\varepsilon_1-\varepsilon_3, \varepsilon_2-\varepsilon_3\},\cr
&\Phi_1^+=\{\delta, \delta\pm \varepsilon_i\mid i=1,2,3\}.
\end{align}

 There is a bilinear form on $\hhh^*$ given by
 $$ (\delta,\delta)=2, (\varepsilon_{i},\varepsilon_{i})=-2, (\varepsilon_{i},\varepsilon_{j})=1, (\delta,\varepsilon_{i})=0, i,j =1,2,3, i\neq j.$$

 The extended standard fundamental system for $G(3)$ is $$\widetilde \Pi=\{-2\delta\}\cup \{\delta-\varepsilon_1,   \varepsilon_2-\varepsilon_3,  -\varepsilon_2\}.$$
 }

All fundamental systems can be classified into following $8$ cases:
\begin{align*}
	& \widetilde\Pi^{'}=\{s(\delta+\varepsilon_{i}),s(\varepsilon_{k}-\varepsilon_{j}),s\varepsilon_{j}\} \cr
	& \widetilde\Pi^{'}=\{s(\delta-\varepsilon_{i}),s(\varepsilon_{j}-\varepsilon_{k}),-s\varepsilon_{j}\} \cr
	& \widetilde\Pi^{'}=\{s(\delta+\varepsilon_{i}),s(-\delta+\varepsilon_{j}),s(\varepsilon_{k}-\varepsilon_{i})\} \cr
	& \widetilde\Pi^{'}=\{s(\delta-\varepsilon_{i}),s(-\delta-\varepsilon_{j}),s(\varepsilon_{i}-\varepsilon_{k})\} \cr
	& \widetilde\Pi^{'}=\{s(\delta+\varepsilon_{i}),s(-\delta-\varepsilon_{j}),-s\varepsilon_{i}\} \cr
	& \widetilde\Pi^{'}=\{s(\delta-\varepsilon_{i}),s(-\delta+\varepsilon_{j}),s\varepsilon_{i}\} \cr
	& \widetilde\Pi^{'}=\{s\delta,s(-\delta+\varepsilon_{i}),s(\varepsilon_{j}-\varepsilon_{i})\} \cr
	& \widetilde\Pi^{'}=\{s\delta,s(-\delta-\varepsilon_{i}),s(\varepsilon_{i}-\varepsilon_{j})\},
\end{align*}
where $s\in\{1,-1\}$, $\{i,j,k\}=\{1,2,3\}$.

\subsection{Coxeter graph}

\begin{theorem}\label{thm: excep g}
	The Coxeter graph of $\wsc(G(3))$  is	
	
	
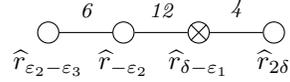
\begin{center}
		\begin{tikzpicture}
			
			\node[circle,draw,minimum size=8pt,inner sep=0pt,fill=white] (r1) at (-4,1) {};

			\node[circle,draw,minimum size=8pt,inner sep=0pt,fill=white] (r2) at (-3,1) {};
			
			\node[circle,draw,minimum size=8pt,inner sep=0pt,fill=white] (r3) at (-2,1) {};
			\draw (r3.135) -- (r3.315);
			\draw(r3.225) -- (r3.45);
			
			\node[circle,draw,minimum size=8pt,inner sep=0pt,fill=white] (r4) at (-1,1) {};

			\draw (r1) -- node[above=2pt] {\footnotesize 6} (r2);
			\draw (r2) -- node[above=2pt] {\footnotesize 12}(r3);
			\draw (r3) --  node[above=2pt] {\footnotesize 4}(r4);

			\node[below=3pt] at (r1) {$ \widehat r_{\varepsilon_2-\varepsilon_3}$};
			\node[below=3pt] at (r2) {$\widehat r_{-\varepsilon_2}$};
			\node[below=3pt] at (r3) {$\widehat r_{\delta-\varepsilon_{1}}$};
			\node[below=3pt] at (r4) {$\widehat r_{2\delta}$};

		\end{tikzpicture}
	\end{center}

\end{theorem}

{
Note that the Coxeter graph of  the subsystem $(\scrw, \spo)$ is already known. It suffices for us to compute $m_{xy}$ involving $\widehat r_{\delta - \varepsilon_1}$. So the orders of $ \widehat r_{\varepsilon_2-\varepsilon_3}\cdot\widehat r_{\delta - \varepsilon_1}$, $\widehat r_{-\varepsilon_2}\cdot\widehat r_{\delta - \varepsilon_1}$ and $ \widehat r_{2\delta}\cdot\widehat r_{\delta - \varepsilon_1}$ should be given.
	
	For	$\widehat r_{-\varepsilon_2}\cdot\widehat r_{\delta - \varepsilon_1}$, we need to calculate the actions of $\widehat r_{-\varepsilon_2}\cdot\widehat r_{\delta - \varepsilon_1}$ on all the  fundamental systems, and get the corresponding orders of them. We also discuss in different cases.

(i) $\widetilde\Pi^{'}$ contains $\pm(\delta-\varepsilon_{1})$ and $\pm(\varepsilon_{2})$, then $\widetilde\Pi^{'}=\{s(\delta-\varepsilon_{1}),s(-\varepsilon_{2}),s(\varepsilon_{2}-\varepsilon_{3})\}$ or $\{s(\delta-\varepsilon_{1}),s(-\varepsilon_{2}),s(-\delta+\varepsilon_{2})\}$, $s\in\{\pm1\}$. The order of both cases is $6$ and we just give the calculation of $\widetilde\Pi^{'}=\{s(\delta-\varepsilon_{1}),s(-\varepsilon_{2}),s(\varepsilon_{2}-\varepsilon_{3})\}$.
\begin{align*}
	&\{s(\delta-\varepsilon_{1}),s(-\varepsilon_{2}),s(\varepsilon_{2}-\varepsilon_{3})\}
	\overset{\widehat r_{\delta-\varepsilon_{1}}}{\longrightarrow}
	\{s(-\delta+\varepsilon_{1}), s(\delta+\varepsilon_{3}),s(\varepsilon_{2}-\varepsilon_{3}) \}
	\overset{\widehat r_{-\varepsilon_{2}}}{\longrightarrow}\cr
&	\{s(-\delta-\varepsilon_{3}), s(\delta-\varepsilon_{1}),s(\varepsilon_{1}-\varepsilon_{2}) \}
	\overset{\widehat r_{\delta-\varepsilon_{1}}}{\longrightarrow}
	\{s\varepsilon_{2}, s(-\delta+\varepsilon_{1}),s(\delta-\varepsilon_{2}) \}\overset{\widehat r_{-\varepsilon_{2}}}{\longrightarrow}\cr
&\{s(-\varepsilon_{2}), s(-\delta-\varepsilon_{3}),s(\delta+\varepsilon_{2})\}
	\overset{\widehat r_{\delta-\varepsilon_{1}}}{\longrightarrow}
	\{s(-\varepsilon_{2}), s(-\delta-\varepsilon_{3}),s(\delta+\varepsilon_{2})\}
	\overset{\widehat r_{-\varepsilon_{2}}}{\longrightarrow}\cr
&\{s\varepsilon_{2}, s(-\delta+\varepsilon_{1}),s(\delta-\varepsilon_{2})\}
	\overset{\widehat r_{\delta-\varepsilon_{1}}}{\longrightarrow}
	\{s(-\delta-\varepsilon_{3}), s(\delta-\varepsilon_{1}),s(\varepsilon_{1}-\varepsilon_{2}) \}
	\overset{\widehat r_{-\varepsilon_{2}}}{\longrightarrow}\cr
&	\{s(-\delta+\varepsilon_{1}), s(\delta+\varepsilon_{3}),s(\varepsilon_{2}-\varepsilon_{3}) \}\overset{\widehat r_{\delta-\varepsilon_{1}}}{\longrightarrow}\{s(\delta-\varepsilon_{1}), s(-\varepsilon_{2}),s(\varepsilon_{2}-\varepsilon_{3}) \}\overset{\widehat r_{-\varepsilon_{2}}}{\longrightarrow}\cr
&	\{s(\delta+\varepsilon_{3}), s\varepsilon_{2},s(\varepsilon_{1}-\varepsilon_{2}) \}\overset{\widehat r_{\delta-\varepsilon_{1}}}{\longrightarrow}\{s(\delta+\varepsilon_{3}), s\varepsilon_{2},s(\varepsilon_{1}-\varepsilon_{2}) \}\overset{\widehat r_{-\varepsilon_{2}}}{\longrightarrow}\cr
&\{s(\delta-\varepsilon_{1}),s(-\varepsilon_{2}),s(\varepsilon_{2}-\varepsilon_{3})\}.
\end{align*}
 So the order of the action  $\widehat r_{\varepsilon_1+\varepsilon_3}\cdot\widehat r_{\delta - \varepsilon_1}$ on $\widetilde\Pi^{'}=\{s(\delta-\varepsilon_{1}),s(-\varepsilon_{2}),s(\varepsilon_{2}-\varepsilon_{3})\}$  is $6$.

(ii) $\widetilde\Pi^{'}$ contains $\pm(\delta-\varepsilon_{1})$ and $\pm(\delta+\varepsilon_{3})$, then $\widetilde\Pi^{'}=\{s(\delta-\varepsilon_{1}),s(-\delta-\varepsilon_{3}),s(\varepsilon_{3}-\varepsilon_{2})\}$ or $\{s(\delta-\varepsilon_{1}),s(-\delta-\varepsilon_{3}),s(\varepsilon_{1}-\varepsilon_{2})\}$, $s\in \{\pm1\}$. The order of both cases is $6$ and we just give the calculation of $\widetilde\Pi^{'}=\{s(\delta-\varepsilon_{1}),s(-\delta-\varepsilon_{3}),s(\varepsilon_{3}-\varepsilon_{2})\}$.

\begin{align*}
	&\{s(\delta-\varepsilon_{1}),s(-\delta-\varepsilon_{3}),s(\varepsilon_{3}-\varepsilon_{2})\}
	\overset{\widehat r_{\delta-\varepsilon_{1}}}{\longrightarrow}
	\{s(-\delta+\varepsilon_{1}), s\varepsilon_{2},s(\varepsilon_{3}-\varepsilon_{2}) \}
	\overset{\widehat r_{-\varepsilon_{2}}}{\longrightarrow}\cr
	&	\{s(-\delta-\varepsilon_{3}), s(-\varepsilon_{2}),s(\varepsilon_{2}-\varepsilon_{1}) \}
	\overset{\widehat r_{\delta-\varepsilon_{1}}}{\longrightarrow}
	\{s(-\delta-\varepsilon_{3}), s(-\varepsilon_{2}),s(\varepsilon_{2}-\varepsilon_{1}) \}\overset{\widehat r_{-\varepsilon_{2}}}{\longrightarrow}\cr
	&\{s(-\delta+\varepsilon_{1}), s\varepsilon_{2},s(\varepsilon_{3}-\varepsilon_{2})\}
	\overset{\widehat r_{\delta-\varepsilon_{1}}}{\longrightarrow}
	\{s(\delta-\varepsilon_{1}), s(-\delta-\varepsilon_{3}),s(\varepsilon_{3}-\varepsilon_{2})\}
	\overset{\widehat r_{-\varepsilon_{2}}}{\longrightarrow}\cr
	&\{s(\delta+\varepsilon_{3}), s(-\delta+\varepsilon_{1}),s(\varepsilon_{2}-\varepsilon_{1})\}
	\overset{\widehat r_{\delta-\varepsilon_{1}}}{\longrightarrow}
	\{s(-\varepsilon_{2}), s(\delta-\varepsilon_{1}),s(-\delta+\varepsilon_{2}) \}
	\overset{\widehat r_{-\varepsilon_{2}}}{\longrightarrow}\cr
	&	\{s\varepsilon_{2}, s(\delta+\varepsilon_{3}),s(-\delta-\varepsilon_{2}) \}\overset{\widehat r_{\delta-\varepsilon_{1}}}{\longrightarrow}\{s\varepsilon_{2}, s(\delta+\varepsilon_{3}),s(-\delta-\varepsilon_{2}) \}\overset{\widehat r_{-\varepsilon_{2}}}{\longrightarrow}\cr
	&	\{s(-\varepsilon_{2}), s(\delta-\varepsilon_{1}),s(-\delta+\varepsilon_{2}) \}\overset{\widehat r_{\delta-\varepsilon_{1}}}{\longrightarrow}\{s(\delta+\varepsilon_{3}), s(-\delta+\varepsilon_{1}),s(\varepsilon_{2}-\varepsilon_{1}) \}\overset{\widehat r_{-\varepsilon_{2}}}{\longrightarrow}\cr
	&\{s(\delta-\varepsilon_{1}),s(-\delta-\varepsilon_{3}),s(\varepsilon_{3}-\varepsilon_{2})\}.
\end{align*}

(iii) $\widetilde\Pi^{'}$ contains $\pm(\delta+\varepsilon_{3})$ and $\pm(\varepsilon_{2})$, then $\widetilde\Pi^{'}=\{s(\delta+\varepsilon_{3}),s(\varepsilon_{1}-\varepsilon_{2}),s\varepsilon_{2}\}$ or $\{s(\delta+\varepsilon_{3}),s\varepsilon_{2},s(-\delta-\varepsilon_{2})\}$, $s\in\{\pm1\}$. Similar to the calculation of action on $\widetilde\Pi^{'}=\{s(\delta-\varepsilon_{1}),s(-\varepsilon_{2}),s(\varepsilon_{2}-\varepsilon_{3})\}$, the order in this case is $6$.

(iv) $\widetilde\Pi^{'}=\{s(\delta-\varepsilon_{1}),\alpha,\beta\}$, where $\alpha,\beta\ne \pm(\varepsilon_{2}),\pm(\delta+\varepsilon_{3})$. By observing all fundamental systems, one can see that $\widehat r_{-\varepsilon_{2}}(\alpha),\widehat r_{-\varepsilon_{2}}(\beta),\widehat r_{-\varepsilon_{2}}\cdot \widehat r_{\delta-\varepsilon_{1}}(\alpha),\widehat r_{-\varepsilon_{2}}\cdot \widehat r_{\delta-\varepsilon_{1}}(\beta)\ne \pm (\delta-\varepsilon_{1})$. Then it follows that
\begin{align*}
	&\{s(\delta-\varepsilon_{1}),\alpha,\beta\}
	\overset{\widehat r_{\delta-\varepsilon_{1}}}{\longrightarrow}
	\{s(-\delta+\varepsilon_{1}), \widehat r_{\delta-\varepsilon_{1}}(\alpha),\widehat r_{\delta-\varepsilon_{1}}(\beta) \}
	\overset{\widehat r_{-\varepsilon_{2}}}{\longrightarrow}\cr
	&	\{s(-\delta-\varepsilon_{3}), \widehat r_{-\varepsilon_{2}}\cdot\widehat r_{\delta-\varepsilon_{1}}(\alpha),\widehat r_{-\varepsilon_{2}}\cdot\widehat r_{\delta-\varepsilon_{1}}(\beta) \}
	\overset{\widehat r_{\delta-\varepsilon_{1}}}{\longrightarrow}
	\{s(-\delta-\varepsilon_{3}), \widehat r_{-\varepsilon_{2}}\cdot\widehat r_{\delta-\varepsilon_{1}}(\alpha),\widehat r_{-\varepsilon_{2}}\cdot\widehat r_{\delta-\varepsilon_{1}}(\beta) \}\overset{\widehat r_{-\varepsilon_{2}}}{\longrightarrow}\cr
	&\{s(-\delta+\varepsilon_{1}), \widehat r_{\delta-\varepsilon_{1}}(\alpha),\widehat r_{\delta-\varepsilon_{1}}(\beta) \}
	\overset{\widehat r_{\delta-\varepsilon_{1}}}{\longrightarrow}
	\{s(\delta-\varepsilon_{1}), \alpha,\beta\}
	\overset{\widehat r_{-\varepsilon_{2}}}{\longrightarrow}\cr
	&\{s(\delta+\varepsilon_{3}), \widehat r_{-\varepsilon_{2}}(\alpha),\widehat r_{-\varepsilon_{2}}(\beta)\}
	\overset{\widehat r_{\delta-\varepsilon_{1}}}{\longrightarrow}
	\{s(\delta+\varepsilon_{3}), \widehat r_{-\varepsilon_{2}}(\alpha),\widehat r_{-\varepsilon_{2}}(\beta) \}
	\overset{\widehat r_{-\varepsilon_{2}}}{\longrightarrow}\cr
	&\{s(\delta-\varepsilon_{1}),\alpha,\beta\}.
\end{align*}
So the order in this case is $4$.

What's more, we can check the following case for order $4$ and we omit the detailed calculation.

(v) $\widetilde\Pi^{'}=\{s(\delta+\varepsilon_{3}),\alpha,\beta\}$, where $\alpha,\beta\ne \pm(\varepsilon_{2}),\pm(\delta-\varepsilon_{1})$. By observation, one can easily see that $\widehat r_{-\varepsilon_{2}}(\alpha),\widehat r_{-\varepsilon_{2}}(\beta),\widehat r_{-\varepsilon_{2}}\cdot \widehat r_{\delta-\varepsilon_{1}}(\alpha),\widehat r_{-\varepsilon_{2}}\cdot \widehat r_{\delta-\varepsilon_{1}}(\beta)\ne \pm (\delta-\varepsilon_{1})$. The order is also $4$.

(vi) $\widetilde\Pi^{'}$ contains no $\pm(\delta+\varepsilon_{3}),\pm(\delta-\varepsilon_{1})$, it is easy to see that the odd reflection acts trivially. Thus in this case, the order is $2$.

 In summary, the order of $\widehat r_{-\varepsilon_2}\cdot\widehat r_{\delta - \varepsilon_1}$ is the least common multiple of $4$ and $6$, which is  equal to $12$.

Next we need to calculate the order of $\widehat r_{\delta}\cdot \widehat r_{\delta-\epsilon_{1}}$. First we notice that if fundamental system $\widetilde\Pi^{'}$ contains no elements $\pm(\delta-\varepsilon_{1})$ and $\pm (\delta+\varepsilon_{1})$, then the order of action $\widehat r_{2\delta}\cdot \widehat r_{\delta-\epsilon_{1}}$ on $\widetilde\Pi^{'}$ is $2$. Now we will consider $2$ cases left.

(i) $\widetilde\Pi^{'}=\{s(\delta-\varepsilon_{1}),\alpha,\beta\}$, then it follows that $\widehat r_{\delta-\epsilon_{1}}(\alpha),\widehat r_{\delta-\epsilon_{1}}(\beta),\widehat r_{\delta}\cdot \widehat r_{\delta-\epsilon_{1}}(\alpha),\widehat r_{\delta}\cdot \widehat r_{\delta-\epsilon_{1}}(\beta)\ne \pm(\delta-\varepsilon_{1})$.
\begin{align*}
	&\{s(\delta-\varepsilon_{1}),\alpha,\beta\}
	\overset{\widehat r_{\delta-\varepsilon_{1}}}{\longrightarrow}
	\{s(-\delta+\varepsilon_{1}), \widehat r_{\delta-\varepsilon_{1}}(\alpha),\widehat r_{\delta-\varepsilon_{1}}(\beta) \}
	\overset{\widehat r_{\delta}}{\longrightarrow}\cr
	&	\{s(\delta+\varepsilon_{1}), \widehat r_{\delta}\cdot\widehat r_{\delta-\varepsilon_{1}}(\alpha),\widehat r_{\delta}\cdot\widehat r_{\delta-\varepsilon_{1}}(\beta) \}
	\overset{\widehat r_{\delta-\varepsilon_{1}}}{\longrightarrow}
	\{s(\delta+\varepsilon_{1}), \widehat r_{\delta}\cdot\widehat r_{\delta-\varepsilon_{1}}(\alpha),\widehat r_{\delta}\cdot\widehat r_{\delta-\varepsilon_{1}}(\beta) \}\overset{\widehat r_{\delta}}{\longrightarrow}\cr
	&\{s(-\delta+\varepsilon_{1}), \widehat r_{\delta-\varepsilon_{1}}(\alpha),\widehat r_{\delta-\varepsilon_{1}}(\beta) \}
	\overset{\widehat r_{\delta-\varepsilon_{1}}}{\longrightarrow}
	\{s(\delta-\varepsilon_{1}), \alpha,\beta\}
	\overset{\widehat r_{\delta}}{\longrightarrow}\cr
	&\{s(-\delta-\varepsilon_{1}), \widehat r_{\delta}(\alpha),\widehat r_{\delta}(\beta)\}
	\overset{\widehat r_{\delta-\varepsilon_{1}}}{\longrightarrow}
	\{s(-\delta-\varepsilon_{1}), \widehat r_{\delta}(\alpha),\widehat r_{\delta}(\beta)\}
	\overset{\widehat r_{\delta}}{\longrightarrow}\cr
	&\{s(\delta-\varepsilon_{1}),\alpha,\beta\}.
\end{align*}
It follows that the order of $\widehat r_{\delta}\cdot \widehat r_{\delta-\epsilon_{1}}$ on this fundamental system is $4$.

Similarly, we can get the order of following case and we omit the details.

(ii) $\widetilde\Pi^{'}=\{s(\delta+\varepsilon_{1}),\alpha,\beta\}$, it is easy to check that the order in this case is also $4$.

Thus we conclude that the order of $\widehat r_{\delta}\cdot \widehat r_{\delta-\epsilon_{1}}$ is $4$.

The calculations are completed.

\section*{Statement}
The authors have no conflict of interest to declare that are relevant to this article.

\bibliographystyle{amsplain}

\begin{thebibliography}{20}
	\bibitem{CLS} C.-J. Chen, Y.-Y. Li and B. Shu, {\em Defining sequences for fundamental root systems and Coxeter graphs for super Weyl groups}, 	Journal of Lie Theory (2026),  to appear  (see arXiv:2401.11068 [math.RT]).
	\bibitem{CW} S.-J. Cheng  and W. Wang,  {\em Dualities and representations of Lie superalgebras}, Graduate Studies in Mathematics, 144. American Mathematical Society, Providence, RI, 2012.
	
	\bibitem{CSW} S.-J. Cheng, B. Shu  and W. Wang,  {\em Modular representations of exceptional supergroups}, Math. Z. 291 (2019), 635-659.








	

\bibitem{GP}  M. Geck and G. Pfeiffer, {\em Characters of finite Coxeter groups and Iwahori-Hecke algebras,} LMS Monogr. (N.S.) 21, The
Clarendon Press, Oxford University Press, New York, 2000.


259 (2008),  255-276.

\bibitem{Hum1} J. E. Humphreys, {\em Introduction to Lie algebras and representation theory}, Grad. Texts in Math. 9, Springer-Verlag, New York-Berlin, 1978.
	
	\bibitem{Hum2} J. E. Humphreys J E. {\em Reflection groups and Coxeter groups}, Cambridge Stud. Adv. Math., 29
Cambridge University Press, Cambridge, 1990.
	


	\bibitem{K} V. G. Kac, {\em Lie superalgebras},  Advances in Math. 26 (1977),  8-96.
	
	




\bibitem{LS} Y.-Y. Li and B.  Shu, {\em
Jantzen filtration of Weyl modules for general
linear supergroups}, Forum Math. 35 (2023), 1435-1468.


	\bibitem{PS} L. Pan and B. Shu, {\em Jantzen filtration and strong linkage principle for modular Lie superalgebras}, Forum Math. 30 (2018), 1573-1598.











	
\end{thebibliography}

\end{document}